\documentclass{amsart}
\usepackage{graphicx}
\newtheorem{theorem}{Theorem}

\usepackage{amsmath,amscd}

\begin{document}

\title{On the Non-Existence of CR-Regular Embeddings of $S^5$}

\author{Ali M. Elgindi}

\begin{abstract}
We show that there exists no CR-regular embedding of the 5-sphere $S^5$ into $\mathbb{C}^4$. We also obtain analogous results for embeddings of higher dimensional spheres into complex space.
\end{abstract}

\maketitle
\par\ \par\

\section{Significance Statement}
\par\ \par\
The study of embeddings of real-submanifolds into complex space is a field of study that lies in the intersection of several branches of mathematics, particularly topology and complex analysis. Much study has been done in the situation of embedding real $n$-dimensional manifolds into $n$-dimensional complex space. Here, we discuss embeddings of spheres into complex spaces of different dimensions, particularly the non-existence of regular embeddings of the $5$-sphere into $4$-dimensional complex space.
\par\

\section{Introduction}
\par\ \par\
The h-principle was developed in the early 1970's by M. Gromov who applied it to the problem of embeddings of real manifolds into complex space to give necessary and sufficient conditions for the existence of totally real embeddings (see [5]). In particular, he demonstrated that the only spheres $S^n$ which admit totally real embeddings to $\mathbb{C}^n$ are in the dimensions $n=1,3$ (the case $n=1$ being trivial).

In the 1980's, F. Forstneric extended the work on Gromov and proved that every compact, orientable 3-manifold admits a totally real embedding into $\mathbb{C}^3$ (see [4]). 
\par\ \par\ 
In our work in [1], we demonstrated that every topological type of knot (or link) in $S^3$ can arise as the set of complex tangents to a $\mathcal{C}^n$-embedding $S^3 \hookrightarrow \mathbb{C}^3$ (for any given integer n). In [2], we derived a topological invariant for the local removal of complex tangents to an embedding of a closed oriented 3-manifold into $\mathbb{C}^3$, leaving the embedding unchanged outside a small neighborhood of a chosen set of complex tangents. This led to our work in [3], where we demonstrated that any embedding $S^3 \hookrightarrow \mathbb{C}^3$ can be approximated $\mathcal{C}^0$-close by a totally real embedding.
\par\ \par\
In their paper [7], N. Kasuya and M. Takase generalized our work in [1] to demonstrate that every knot or link in $S^3$ can be exactly assumed as the set of complex tangents to a smooth embedding of $S^3$. In fact, they generalize further to demonstrate that any 1-dimensional submanifold of a closed orientable 3-manifold that is homologically trivial may arise as the set of complex tangents to an embedding $M \hookrightarrow \mathbb{C}^3$. Their method of proof is based on the theory of stable maps in topology and Saeki's Theorem.
\par\ \par\
In more general situations, in particular for embeddings $f: M^n \hookrightarrow \mathbb{C}^q$ ($n \neq q$) the h-principle still holds, although there cannot be in the literal sense totally real embeddings if $n>q$. In the situation $n>q$, every point $x \in M$ must be complex tangent, and the $\textit{complex dimension}$ of $x$ is defined to be: 
\par\
$dim(x) = dim_{\mathbb{C}} (f_*(T_x M) \bigcap J f_* (T_x M))$ where $J$ is the complex structure. Note that $dim(x) \geq n-q$, by elementary linear algebra. 
\par\
If $dim(x) = n-q$, we say that $x$ is a $\textit{CR-regular}$ point of the embedding. If $dim(x) > n-q$, we say $x$ is called $\textit{CR-singular}$.
\par\ \par\
In his paper [9], M. Slapar considered the situation of a closed oriented 4-manifold embedded into $\mathbb{C}^3$. In this situation, complex tangents generically are discrete (and finite), and can each be given sign (plus or minus) by comparing the orientation of the tangent space of the embedded manifold at the complex tangent with the induced orientation of the tangent space as a complex subspace of the tangent space of the complex ambient manifold.  Slapar also demonstrated (in [8]) that a 4-manifold admits a CR-regular embedding into $\mathbb{C}^3$ if and only if the manifold is parallelizable.
\par\ \par\
In this paper, we will focus on embeddings $S^5 \hookrightarrow \mathbb{C}^4$. In this situation, a generic embedding will assume CR-singular points along a knot in $S^5$, all other points being CR-regular. We will show that there exists no CR-regular embedding of $S^5 \hookrightarrow \mathbb{C}^4$, and we will also make further generalizations for embeddings of spheres of higher dimensions.
\par\
\section{The Result and Proof}
\par\ \par\
An embedding of $S^5 \hookrightarrow \mathbb{C}^4$ must have a complex line in the tangent space to each point $x \in S^5$, because of dimensionality reasons. There are two classes of points in such an embedding of $S^5$. In particular a point may have a complex plane in its tangent space, we say such a point is CR-singular. If the tangent plane at the point does not contain a complex plane in its tangent space (only a complex line), we say the point is CR-regular. For a generic embedding  $S^5 \hookrightarrow \mathbb{C}^4$, we will have a link of CR-singular points, and all other points will be CR-regular (this readily follows from the dimension of the relevant spaces and using an application of Sard's theorem).
\par\ \par\
Now, let $\mathbb{G}_{5,8}$ be the Grassmannian of 5-planes in $\mathbb{R}^8 = \mathbb{C}^4$, and consider its subspace $\mathbb{Y} = \{P \in \mathbb{G}_{5,8} | P \bigcap JP  \cong \mathbb{C} \}$ of planes that contain only a complex line. We say a 5-plane $P \in \mathbb{Y}$ is $\textit{almost real}$.
\par\ \par\
Consider the Gauss map of the embedding: $\textit{G}: S^5 \rightarrow \mathbb{G}$. Then the CR-regular points of the embedding are exactly the points whose image under $\textit{G}$ is contained in $\mathbb{Y}$. In particular:
\par\
$\{$CR-regular points$\} = \textit{G}^{-1} (\mathbb{Y}) \subset S^5$
\par\ \par\
A CR-regular embedding $S^5 \hookrightarrow \mathbb{C}^4$ will have that its Gauss map has image only in $\mathbb{Y}$, that is the tangent plane at every point is totally real.
\par\ 
 We will show that no embedding can satisfy this, in particular:
\par\ 
\begin{theorem}:
There exists no CR-regular embedding $S^5 \hookrightarrow \mathbb{C}^4$.
      \end{theorem}
\par\ \par\ 
\textit{\textbf{\underline{Proof}}}:-
\par\ 
As above, we denote $\mathbb{Y} \subset \mathbb{G}_{5,8}$ to be the space of almost real 5-planes in $\mathbb{R}^8 \cong \mathbb{C}^4$, i.e. a 5-plane $\mathcal{P} \in \mathbb{Y}$ if and only if $\mathcal{P} \cap J(\mathcal{P}) \cong \mathbb{C}$.
\par\ \par\
Consider such a 5-plane $\mathcal{P} \in \mathbb{Y}$, and let $L = \mathcal{P} \cap J(\mathcal{P})$ be the (unique) complex line it contains. Further, denote the complement of $L$ in $\mathcal{P}$ by $\mathcal{M}$, so now: $\mathcal{P} = L \oplus \mathcal{M}$. Note that $\mathcal{M} \subset \mathbb{C}^4$ is necessarily a totally real 3-plane (as $\mathcal{P} \in \mathbb{Y}$ contains only a complex line in $L$).
\par\ \par\
Further denote: $\mathcal{N}=\mathcal{P}^\perp$ to be the normal 3-plane to $\mathcal{P} \subset \mathbb{R}^8$. 
\par\ \par\
We note the following from linear algebra:
\par\
$(\mathcal{M} \oplus J\mathcal{M})^\perp = \mathcal{M}^\perp \cap (J\mathcal{M})^\perp = \mathcal{M}^\perp \cap J(\mathcal{M}^\perp) = \mathcal{P} \cap J(\mathcal{P}) = L$.
\par\ \par\
It is then clear that $\mathcal{M} \oplus J\mathcal{M} = L^\perp = \mathcal{M} \oplus \mathcal{N}$. Hence, we have that the normal plane $\mathcal{N} = J(\mathcal{M})$.  
\par\ \par\
Suppose now that we had a CR-regular embedding $F: S^5 \hookrightarrow \mathbb{C}^4$. Note that using this CR-regular embedding we can decompose the tangent bundle of $S^5$ as the Whitney sum of two bundles, namely a 2-plane bundle given by the complex lines at each point using the regular CR-structure and the 3-plane bundle of orthogonal totally real 3-planes at each point. 
\par\
Using notation analogous to our above work (but now using bundles), we may write: $T(S^5) \cong L \oplus \mathcal{M}$. 
\par\ \par\
We note that the normal bundle to any embedding $S^5 \hookrightarrow \mathbb{R}^8$ is trivial (see Massey in [8]). Let $\mathcal{N}$ denote the normal bundle to this embedding, and note by the above work that $\mathcal{M} = J(\mathcal{N})$ as bundles. Now, with a choice of trivialization of the normal bundle we may precompose with the complex structure $J$ to obtain a trivialization of the bundle $\mathcal{M}$.
\par\ \par\
Now, it is well-known that $k$-plane bundles over spheres $S^n$ are classified by the homotopy groups $\pi_{n-1}(GL(k))$ (see Hatcher in [6]). Hence, $2$-plane bundles ($k=2$) over spheres are classified by the homotopy groups of $GL(2) \cong S^1$, and so $2$-plane bundles over spheres of dimension $\geq 3$ are all trivial. Particularly in our case (over $S^5$) the bundle $L$     must be trivial.
\par\ \par\
Therefore, $T(S^5) \cong L \oplus \mathcal{M}$ is equivalent to a sum of trivial bundles. But the Whitney sum of trivial bundles must be trivial! 
\par\ \par\
This gives a contradiction, as the 5-sphere is not parallelizable. So we may conclude that there exists no CR-regular embedding of $S^5 \hookrightarrow \mathbb{C}^4$.
\par\ \par\
$\textbf{\textit{QED}}$
\par\ \par\
We now recall the following result of Kervaire (see Massey in [8]):
\par
The normal bundle to an $n$-sphere embedded in $\mathbb{R}^{n+k}$ with $k>\frac{n+1}{2}$ is necessarily trivial.
\par\ \par\
We now generalize our above work to achieve higher dimensional analogues:
\par\ 
\begin{theorem}:
There exists no CR-regular embedding $S^n \hookrightarrow \mathbb{C}^{n-1}$ for $n \geq 6 , n \neq 7$.
      \end{theorem}
\par\  \par\
\textit{\textbf{\underline{Proof}}}:-
\par\ 
Applying the theorem of Kervaire for $k = n-2$, we find that the normal bundle of an $n$-sphere embedded in $\mathbb{R}^{2n-2}$ is necessarily trivial. Note that we have achieved the analogous result for the case $n=5$ in our first theorem above (which is not covered in the dimension range of this theorem).
\par\ \par\
Suppose there exists a CR-regular embedding $F:S^n \hookrightarrow \mathbb{C}^{n-1}$, and let $x \in S^n$. We would then get that the holomorphic tangent space $T_x \cap JT_x = L_x$ is a complex line. Analogous to our above work, we may consider the complement $\mathcal{M}_x = L^\perp \subset T_x$ which is necessarily a totally real $(n-2)$-plane in $\mathbb{C}^{n-1}$. These will form bundles $L$ and $\mathcal{M}$ over $S^5$. Let $\mathcal{N}$ denote the normal bundle to this embedding.
\par\ \par\
Again, it is easy to see that $\mathcal{M} = J (\mathcal{N})$ and so the trivialization of the normal bundle will give a trivialization of $\mathcal{M}$. Also, the real 2-plane bundle $L$ is necessarily trivial as it is a 2-plane bundle over $S^n$, $n > 3$. 
\par\ \par\
Now, as we may now decompose $T (S^n) = L \oplus \mathcal{M}$ as a sum of two trivial bundles, the tangent bundle must also be trivial. But $S^n$ is not parallelizable for $n \neq 3, 7$!
\par\ \par\
This gives a contradiction to our assumption of such a CR-regular embedding, and thus proves that there exists no CR-regular embedding $S^n \hookrightarrow \mathbb{C}^{n-1}$ in the dimensions $n \geq 6 , n \neq 7$.
\par\ \par\
$\textbf{\textit{QED}}$
\par\ \par\


\begin{thebibliography}{9}

\bibitem{[1]} A.M. Elgindi, "On the Topological Structure of Complex Tangencies to Embeddings of $S^3$ into $\mathbb{C}^3$," New York J. of Math. Vol. 18 (2012), 295-313.

\bibitem{[2]} A.M. Elgindi "A topological obstruction to the removal of a degenerate complex tangent and some related homotopy and homology groups," International Journal of Mathematics Vol. 26 (2015).

\bibitem{[3]} A.M. Elgindi, "Totally real perturbations and nondegenerate embeddings of $S^3$," New York J. of Math. Volume 21 (2015) 1283-1293

\bibitem{[4]} F. Forstneric, "On totally real embeddings into $\mathbb{C}^n$," Expositiones Mathematicae, 4 (1986), pp. 243-255.

\bibitem{[5]} M.L. Gromov, "Convex Integration of Differential Relations," 1973 Math. USSR Izv. 7.

\bibitem{[6]} A. Hatcher, "Algebraic Topology," Cambridge University Press (2002).

\bibitem{[7]} N. Kasuya and M. Takase, "Knots and links of complex tangents." arXiv preprint 1606.03704 (2016).

\bibitem{[8]} W.S. Massey, "On the Normal Bundle of a Sphere Imbedded in Euclidean Space," Proc. of the Amer. Math. Soc., Vol. 10, No. 6 (Dec., 1959), pp. 959-964.

\bibitem{[9]} M. Slapar, "Cancelling complex points in codimension two," Bull. Aust. Math. Soc. 88 (2013), no. 1, 64-69.


\end{thebibliography}
\end{document}